\documentclass[manmat, numbook, envcountsame, draft]{svjour}
\usepackage{latexsym, amssymb, theorem, amsmath}

\def\CM{Cohen--Macaulay }
\def\Cl{\mathop{\mathrm{Cl}}}

\def\sep{^{+{\mathrm{sep}}}}
\def\gr{^{+{\mathrm{GR}}}}
 
\begin{document}

\title{Separable integral extensions and plus closure}

\author{Anurag K. Singh}

\institute{Anurag K. Singh \at Department of Mathematics, University of
Illinois at Urbana-Champaign, 1409 W. Green Street, Urbana, IL 61801, USA \\
e-mail: {\tt singh6@math.uiuc.edu}}

\date{Received: 28 May 1998 / Revised version: 30 November 1998}


\maketitle

\begin{abstract}
We show that an excellent local domain of characteristic $p$ has a separable
big Cohen--Macaulay algebra. In the course of our work we prove that an element
which is in the Frobenius closure of an ideal can be forced into the expansion
of the ideal to a module--finite separable extension ring. 
\end{abstract}

\section{Introduction}
\label{intro}

Let $R$ be an excellent domain of characteristic $p$, and let $R^+$ denote the 
integral closure of $R$ in an algebraic closure of its fraction field. A
celebrated theorem of M.~Hochster and C.~Huneke states that $R^+$ is a big \CM
algebra for $R$, see \cite{HHbigcm}. The {\it plus closure}\/ $I^+$ of an ideal
$I$ of $R$ is defined as $I^+ = IR^+ \cap R$, and has close connections with
tight closure theory: the containment $I^+ \subseteq I^*$ is easily verified,
and in \cite{Sminv} K.~E.~Smith showed that $I^+ = I^*$ if the ideal $I$ is 
generated by part of a system of parameters for $R$. In this paper we establish
the somewhat surprising result that for any element $z \in I^+$ there exists an
integral domain $S$, which is a {\it separable} module--finite extension of
$R$, such that $z \in IS$. We use this idea to obtain a separable big \CM
algebra $R\sep$ for any excellent local domain $(R,m)$ of characteristic $p$.

For related work on $R^+$ and plus closure see \cite{Aberbach} and \cite{AH}.
Our references for the theory of tight closure are \cite{HHjams},
\cite{HHbasec} and \cite{HHjalg}.

In Section \ref{plus} we present an explicit computation of plus closure,
specifically we show that $xyz \in (x^2, y^2, z^2)^+$ in the cubic hypersurface
$R=K[X,Y,Z]/(X^3+Y^3+Z^3)$. This strengthens a result obtained in \cite{cubic}
where it was established that $xyz \in (x^2, y^2, z^2)^*$.

\section{Preliminaries}
\label{prelim}

Let $R$ be a Noetherian ring of characteristic $p > 0$. The letter $e$ shall 
denote a variable nonnegative integer, and $q$ shall denote the corresponding
power of the characteristic, i.e., $q=p^e$. For an ideal $I=(x_1,\dots, x_n)
\subseteq R$, let $I^{[q]} = (x_1^q,\dots, x_n^q)$. For a reduced ring $R$ of
characteristic $p > 0$, $R^{1/q}$ shall denote the ring obtained by adjoining
all $q\,$th roots of elements of $R$.

For an element $x$ of $R$ we say that $x \in I^F$, the {\it Frobenius
closure}\/ of  $I$, if there exists an integer $q=p^e$ such that $x^q \in
I^{[q]}$. The ring $R$ is said to be {\it F--pure }\/ if the Frobenius
homomorphism is pure, i.e., if $F: M \to F(M)$ is injective for all
$R$--modules $M$.

In specific examples which are homomorphic images of polynomial rings, we shall
use lower case letters to denote the images of the corresponding  variables.

\section{Separable integral extensions}
\label{sepint}

Our main result regarding Frobenius closure and separable extensions is the
following theorem:

\begin{theorem} 
Let $R$ be a excellent domain of characteristic $p>0$, $I\subseteq R$ an ideal, 
and $z\in R$ an element such that $z \in I^F$. Then there exists an integral 
domain $S$, which is a module--finite separable extension of $R$, such that 
$z \in IS$. 
\label{main}
\end{theorem}

\begin{proof}
Since $z \in I^F$, there exists a positive integer $q=p^e$, nonzero elements 
$x_0, \dots, x_n \in I$, and $a_0, \dots, a_n \in R$ such that
$$
z^q = \sum_{i=0}^{n} a_i x_i^q.
$$ 
For $1 \le i \le n$, consider the equations
$$
U_i^q+U_ix_0^q-a_i = 0.
$$
These are monic separable equations in the variables $U_i$, and therefore have
solutions $u_i$ in a separable field extension of the fraction field of $R$. 
Let $S$ be the integral closure of the ring $R[u_1, \dots, u_n]$ in its field
of fractions. We shall show that $z \in IS$. Let
$$
u_0 = (z- \sum_{i=1}^{n} x_i u_i) / x_0,
$$ 
which is an element of the fraction field of $S$. Taking the $q$\,th power of
this element, we have
\begin{align*}
u_0^q \ &= \ ( z^q-\sum_{i=1}^{n}x_i^q u_i^q) / x_0^q \ = \ 
(\sum_{i=0}^{n} a_i x_i^q - \sum_{i=1}^{n} x_i^q u_i^q) / x_0^q \\
\ &= \ a_0 + \sum_{i=1}^{n}(a_i-u_i^q)x_i^q/x_0^q \ = \ 
a_0 + \sum_{i=1}^{n} u_i x_i^q.
\end{align*}
Consequently $u_0$ is integral over $S$, but since $S$ is a normal domain, 
we then have $u_0 \in S$. This implies that $z=\sum_{i=0}^{n} x_i u_i \in IS$.
\qed 
\end{proof}

\begin{remark} 
Let $GL_n({\mathbb F}_q)$ be the general linear group over the finite field 
${\mathbb F}_q$. In \cite{failure} the author constructed examples to  show
that the ring of invariants for the natural action of a subgroup  $G$ of
$GL_n({\mathbb F}_q)$ on the polynomial ring ${\mathbb F}_q[X_1, \dots, X_n]$
need not be F--pure. Specifically, consider the natural action of the
symplectic group  $Sp_4({\mathbb F}_q)$ on the polynomial ring ${\mathbb
F}_q[X_1,X_2,X_3,X_4]$.  Then the ring of invariants is isomorphic to the
hypersurface
$$
R = {\mathbb F}_q[X,Y,Z,A,B]/(Z^q-AX^q-BY^q)
$$ 
which is not F--pure since $z \in I=(x,y)^F$, more precisely the element $z$ is
forced into the expanded ideal $IR^{1/q}$ in the purely inseparable extension
$R^{1/q}$. However $z$ is also forced into the expansion of $I$ to the linearly
disjoint separable extension ${\mathbb F}_q[X_1,X_2,X_3,X_4]$. These examples
provided the basic case of the equational construction used in the proof of
Theorem \ref{main} above.    
\end{remark}

\begin{definition}
Let $R$ be an integral domain with fraction field $\mathcal K$. Let 
$\overline{\mathcal K}$ be an algebraic closure of the field $\mathcal K$, and
$\mathcal L$ be the maximal separable field extension of $\mathcal K$ in
$\overline{\mathcal K}$. Then $R\sep$ shall denote the integral closure of $R$ 
in $\mathcal L$.   
\end{definition}

With this notation, Theorem \ref{main} gives us the following corollary:

\begin{corollary} 
Let $R$ be a excellent integral domain of characteristic $p$. If $I$ is an
ideal of $R$ and $z$ an element of $R$ such that $z \in I^+$, then there exists 
a module--finite separable extension domain $S$ such that $z \in IS$.
Consequently for all ideals $I$ of $R$ we have
$$
I R^+ \cap R = I R\sep \cap R.
$$
\label{closure}
\end{corollary}

\begin{proof}
Since $z \in I^+$ there exists a integral domain $R_1$, module--finite over
$R$, such that $z \in IR_1$. If $R_2$ denotes the largest separable extension
of $R$ contained in  $R_1$, we have $z \in (IR_2)^F$. By Theorem \ref{main}
there exists a separable module--finite extension $S$ of $R_2$ such that $z \in
IS$. 
\qed 
\end{proof}

Let $R = \oplus_{i \in \mathbb N}R_i$ be an $\mathbb N$--graded integral domain
over a field  $K=R_0$ of characteristic $p$. An element $z \in R^+$ is said to
be {\it homogeneous}\/ if it satisfies an equation of integral dependence over
$R$ of the form 
$$ 
z^n + a_1z^{n-1} + \dots + a_n =0 
$$  
where $a_i \in R_{id}$ for some $d \in \mathbb Q$, for all $1 \le i \le n$.
(The element $z$ can then be assigned weight $d$, and satisfies a homogeneous
equation of integral dependence over $R$.) Let $R\gr$ denote the $\mathbb
Q$--graded sub\-algebra of $R^+$ which is generated by homogeneous elements $z
\in R^+$.  In this setting, Hochster and Huneke show that $R\gr$ is a graded
big \CM algebra for $R$, i.e., that every homogeneous system of parameters for
$R$ is a regular sequence on $R\gr$, see \cite{HHbigcm}.

The following example shows that for $R$ as above, a homogeneous ideal $I$ of
$R$, and a homogeneous element $\xi \in I^+$, there need not exist a graded
separable module--finite extension $S$ with $\xi \in IS$.

\begin{example}
Let $R=K[X, Y, Z]/(X^3+Y^3+Z^3)$ where $K$ is an algebraically closed field of 
characteristic $2$. Note that $z^2 \in (x,y)^F$ since $z^4 = zx^3+zy^3$. We
claim $R$ has no module--finite graded separable extension $S$ with $z^2 \in
IS$. If $S$ were such an extension, there exist homogeneous elements $u,v \in S$ 
(of weight $1$) with $z^2=ux+vy$. In $R\gr$ we then have
$$
z^2=ux+vy = x \sqrt{xz} + y \sqrt{yz},
$$
and since $R\gr$ is a graded big \CM algebra for $R$, the relation
$x(u+\sqrt{xz})=y(v+\sqrt{yz})$ must be  trivial in $R\gr$. Hence there exists
$c \in R\gr $ such that $u+\sqrt{xz}=cy$ and $v+\sqrt{yz}=cx$, but then $c$
must have weight $0$, and so satisfies an equation of integral dependence over
the algebraically closed field $K$. Consequently $c \in K$, and therefore
$\sqrt{xz}$ and $\sqrt{yz}$ are elements of a separable extension of $R$, a
contradiction.   
\end{example}

\begin{remark}
Let $R$ be a ring of characteristic $p$ which is not F--pure, i.e., such that
there exists $z \in R$ with $z \in I^F - I$. By Theorem \ref{main} there exists
a module--finite separable extension $S$ such that $z \in IS$. In general we 
cannot expect $S$ to be F--regular or even F--rational: if $R$ is the 
coordinate ring of an elliptic curve, a normal module--finite extension domain
$S$ will continue to have elements of infinite order in the divisor class
group  $\Cl(S)$, and consequently the two dimensional ring $S$ cannot be
F--rational by a result of J.~Lipman, \cite{Li}. However there do exist 
examples where the separable extension is F--regular, moreover is a polynomial 
ring.  
\end{remark}

\begin{example}
Let $R = K[X,Y,Z,A,B]/(Z^q-AX^q-BY^q)$ where $K$ is a field of characteristic
$p$ and $q=p^e$. Then $z \in (x,y)^F$, and we construct a separable
module--finite extension $S$ such that $z \in (x,y)S$. Let $u$ be a root of the 
separable equation 
$$
U^q -a +Uy^q = 0
$$
and let $S$ be the normalization of $R[u]$. Then $v = (z-ux)/y$ is easily seen
to be an element of $S$ and we have $S= K[X,Y,u,v]$, which is a polynomial ring. 
\end{example}

We next recall a definition from \cite{Ma}:

\begin{definition}
A Noetherian integral domain $R$ is said to be a {\it splinter}\/ if it is a
direct summand, as an $R$--module, of every module--finite extension ring of $R$. 
\end{definition}

In the case that $R$ contains the field of rational numbers, it is easily seen
that $R$ is a splinter if and only if it is a normal ring, but the notion is
more subtle for rings of characteristic $p$. F--regular rings are always
splinter and the converse is known to hold for $\mathbb Q$--Gorenstein rings,
see \cite{splinter}. Corollary \ref{closure} above gives a new characterization
of splinter rings of characteristic $p$:

\begin{corollary} 
Let $R$ be an excellent integral domain on characteristic $p$. Then $R$ is a
splinter if and only if it is a direct summand of every module--finite
separable extension domain.  
\end{corollary}

\begin{proof}
Excellent integral domains are approximately Gorenstein, and so an inclusion $R
\to S$  splits if and only if $IS \cap R = I$ for all ideals $I$ of $R$, see 
\cite{approx-gor}. The result now follows from Corollary \ref{closure}.  
\qed
\end{proof}

We also obtain the following theorem which is separable analogue of the main
result of \cite{HHbigcm}.

\begin{theorem} 
Let $(R,m)$ be an excellent local domain of characteristic $p$. Then every
sequence of elements which is part of a system of parameters for $R$
is a regular sequence on $R\sep$. Consequently $R\sep$ is a balanced big
\CM module for $R$. 
\end{theorem}

\begin{proof}
Let $x_1, \dots, x_k$ be part of a system of parameters for $R$. Given a
relation $\sum_{i=1}^k r_i x_i = 0$ with $r_i \in R\sep$ we may replace $R$ by
a module--finite separable extension and (after a change of notation) assume
that $r_i \in R$. Since $R^+$ is a big \CM algebra for $R$, \cite[Theorem
1.1]{HHbigcm}, we have $r_k \in (x_1, \dots, x_{k-1})R^+ \cap R$, but then 
$r_k \in (x_1, \dots, x_{k-1})R\sep$ by Corollary \ref{closure}.  
\qed
\end{proof}

\section{A computation of plus closure}
\label{plus}

In \cite{cubic} we showed that $xyz \in (x^2, y^2, z^2)^*$ in the cubic
hypersurface $R=K[X,Y,Z]/(X^3+Y^3+Z^3)$ where $K$ is a field of prime
characteristic $p \ne 3$. This question arose in M.~McDermott's study of the
tight closure and plus of various irreducible ideals in $R$, see \cite{moira}.
We furthermore showed that $xyz \in (x^2, y^2, z^2)^F$ whenever $R$ is not
F--pure, i.e., when $p \equiv 2 \mod 3$. In this section we settle the one
unresolved issue by establishing that $xyz \in (x^2, y^2, z^2)^+$ when the
characteristic of the field is $p \equiv 1 \mod 3$.

The following lemma is a basic version of a more powerful equational
criterion, but will suffice for our needs.

\begin{lemma} 
Let $R$ be a integral domain of characteristic $p$ and consider nonzero
elements $z, x_1, \dots, x_k \in R$  which satisfy
$$
z^p \in (x_1^p, \dots, x_k^p) + z((x_1^p, \dots, x_k^p) : (x_1, \dots, x_k)).
$$
Then there exists an integral domain $S$ which is a module--finite extension of 
$R$ such that $z \in (x_1, \dots, x_k)S$.
\label{equational}
\end{lemma}

\begin{proof}
See \cite{Ianthesis}.
\qed 
\end{proof}

We record a determinant identity that we shall find useful. For
integers $n$ and $m$ where $m \ge 1$, we shall use the notation:
$$
\binom{n}{m} = \frac{(n)(n-1) \dotsm (n-m+1)}{(m)(m-1) \dotsm (1)}.
$$

\begin{lemma}
$$
\det\begin{pmatrix}
\binom{n}{a+k} & \binom{n}{a+k+1} & \hdots & \binom{n}{a+2k} \\
\binom{n}{a+k-1} & \binom{n}{a+k} & \hdots & \binom{n}{a+2k-1} \\
\hdotsfor4 \\
\hdotsfor4 \\
\binom{n}{a} & \binom{n}{a+1} & \hdots & \binom{n}{a+k} \\
\end{pmatrix} 
=\frac{\binom{n}{a+k}\binom{n+1}{a+k}\dotsm\binom{n+k}{a+k}}
{\binom{a+k}{a+k} \binom{a+k+1}{a+k}\cdots\binom{a+2k}{a+k}}.
$$
\label{det}
\end{lemma}

\begin{proof}
This is evaluated in \cite[page 682]{muir} as well as \cite{roberts}.
\qed \end{proof}

As a first application of this, we prove the following lemma:

\begin{lemma}
Let $K[A,B]$ be a polynomial ring over a field $K$ of characteristic $p=3k+1$,
where $k$ is a positive integer. Then
$$
(A,B)^{3k} \subseteq I = (A^{2k+1}, \ B^{2k+1}, \ (A+B)^{2k}).
$$
\label{general}
\end{lemma}

\begin{proof}
Note that $I$ contains the following elements:
$$
(A+B)^{2k}A^{k}, \ (A+B)^{2k}A^{k-1}B, \dots, \ (A+B)^{2k}B^{k}.
$$
We consider the binomial expansions of these elements modulo the ideal 
$(A^{2k+1}, B^{2k+1})$. Consequently the following elements are in $I:$
$$
\begin{array}{rrrrrrr}
\binom{2k}{k}A^{2k}B^{k}   &+ &\binom{2k}{k+1}A^{2k-1}B^{k+1} &+ &\cdots
&+ &\binom{2k}{2k}A^{k}B^{2k}, \\
\binom{2k}{k-1}A^{2k}B^{k} &+ &\binom{2k}{k}A^{2k-1}B^{k+1}   &+ &\cdots
&+ &\binom{2k}{2k-1}A^{k}B^{2k}, \\
\hdotsfor7 \\
\hdotsfor7 \\
\binom {2k}{0}A^{2k}B^{k}  &+ &\binom{2k}{1}A^{2k-1}B^{k+1}   &+ &\cdots 
&+ &\binom{2k}{k}A^{k}B^{2k}.
\end{array}
$$
The coefficients of $A^{2k}B^{k}, \ A^{2k-1}B^{k+1}, \dots, \ A^{k}B^{2k}$ 
form the matrix: 
$$
\begin{pmatrix}
\binom{2k}{k}   &\binom{2k}{k+1} &\cdots &\binom{2k}{2k} \\
\binom{2k}{k-1} &\binom{2k}{k}   &\cdots &\binom{2k}{2k-1} \\
\hdotsfor4 \\
\hdotsfor4 \\
\binom {2k}{0}  &\binom{2k}{1}   &\cdots &\binom{2k}{k}
\end{pmatrix}.
$$
To show that all monomials of degree $3k$ in $A$ and $B$ are in $I$, it 
suffices to show that this matrix is invertible. By Lemma \ref{det}, the 
determinant of this matrix is 
$$ 
\frac{\binom{2k}{k}\binom{2k+1}{k}\cdots\binom{3k}{k}}
{\binom{k}{k} \binom{k+1}{k}\cdots\binom{2k}{k}}
$$
which is immediately seen to be invertible since the characteristic of the field 
is $p=3k+1$.
\qed \end{proof}

\begin{lemma} 
Let $R=K[X,Y,Z]/(X^3+Y^3+Z^3)$ where $K$ is a field of prime characteristic 
$p=3k+1$. Then
$$
x^{p-1} y^{2p-2} \in (x^{2p}, y^{2p}, z^{2p}):(x^2, y^2, z^2).
$$
\label{colon}
\end{lemma}

\begin{proof}
We first show that $x^{p+1} y^{2p-2} \in (x^{2p}, y^{2p}, z^{2p})$, i.e., that
$$
x^{3k+2} y^{6k} \in (x^{6k+2}, \ y^{6k+2}, \ z^{6k+2}).
$$ 
This would follow if we establish that 
$$
x^{3k} y^{6k} \in (x^{6k}, \ y^{6k+3}, \ z^{6k+3}).
$$
Setting $A=y^3$, $B=z^3$ and $A+B = -x^3$, we need to show 
$$
(A+B)^k A^{2k} \in I = ((A+B)^{2k}, \ A^{2k+1}, \ B^{2k+1}).
$$
This follows immediately from Lemma \ref{general}.

It remains to check that 
$x^{p-1} y^{2p-2}z^2 \in (x^{2p}, y^{2p}, z^{2p})$. This would follow if we 
show
$$
x^{3k} y^{6k} \in (x^{6k+3}, \ y^{6k+3}, \ z^{6k}).
$$
Setting $A=x^3, B=y^3$ and $A+B = -z^3$, we now need to show 
$$
A^k B^{2k} \in I = (A^{2k+1}, \ B^{2k+1}, \ (A+B)^{2k}),
$$
which again follows from Lemma \ref{general}.
\qed \end{proof}

We are now ready to prove our main result of this section.

\begin{theorem}
Let $R=K[X,Y,Z]/(X^3+Y^3+Z^3)$ where $K$ is a field of prime characteristic 
$p \equiv 1 \mod 3$. Then $xyz \in (x^2, y^2, z^2)^+$. 
\end{theorem}

\begin{proof} 
By Lemma \ref{equational}, it suffices to show
$$
(xyz)^p \in (x^{2p}, y^{2p}, z^{2p}) 
 + xyz ((x^{2p}, y^{2p}, z^{2p}) : (x^2, y^2, z^2)).
$$
By Lemma \ref{colon} we have
$$
(x^{p-1}y^{2p-2}, \ x^{2p-2}y^{p-1} ) \subseteq 
  (x^{2p},y^{2p},z^{2p}):(x^2, y^2, z^2)),
$$
and hence it is enough to show
$$
(xyz)^p \in (x^{2p},\  y^{2p}, \ z^{2p}, \ x^{p}y^{2p-1}z, \ x^{2p-1}y^{p}z),
$$
i.e., that
$$
(xyz)^{3k+1} \in (x^{6k+2}, \ y^{6k+2}, \ z^{6k+2}, 
 \ x^{3k+1}y^{6k+1}z, \ x^{6k+1}y^{3k+1}z).
$$
This follows if we show
$$
(xyz)^{3k} \in (z^{6k+3}, \ x^{3k}y^{6k}, \ x^{6k}y^{3k}).
$$
Setting $A=x^3, B=y^3$ and $A+B = -z^3$, we need to show 
$$
(AB(A+B))^k \in I = ((A+B)^{2k+1}, \ A^k B^{2k}, \ A^{2k} B^k).
$$
Note that $I$ contains the following elements:
$$
(A+B)^{2k+1}A^{k-2}B, \ (A+B)^{2k+1}A^{k-3}B^2, \dots, \ (A+B)^{2k+1}B^{k-1}.
$$
We consider the binomial expansions of these elements modulo the ideal 
$(A^{k}B^{2k}, A^{2k}B^{k})$. This shows that the following are elements of $I:$
$$
\begin{array}{rrrrrrr}
\binom{2k+1}{k}A^{2k-1}B^{k+1}   &+ &\binom{2k+1}{k+1}A^{2k-2}B^{k+2}&+ &\cdots
&+ &\binom{2k+1}{2k-2}A^{k+1}B^{2k-1}, \\
\binom{2k+1}{k-1}A^{2k-1}B^{k+1} &+ &\binom{2k+1}{k}A^{2k-2}B^{k+2}  &+ &\cdots
&+ &\binom{2k+1}{2k-3}A^{k+1}B^{2k-1}, \\
\hdotsfor7 \\
\hdotsfor7 \\
\binom {2k+1}{2}A^{2k-1}B^{k+1}  &+ &\binom{2k+1}{3}A^{2k-2}B^{k+2}  &+ &\cdots 
&+ &\binom{2k+1}{k}A^{k+1}B^{2k-1}.
\end{array}
$$
The coefficients of $A^{2k-1}B^{k+1}, \ A^{2k-2}B^{k+2}, \dots, \ A^{k+1}B^{2k-1}$ 
form the matrix: 
$$
\begin{pmatrix}
\binom{2k+1}{k}   &\binom{2k+1}{k+1} &\cdots  &\binom{2k+1}{2k-2} \\
\binom{2k+1}{k-1} &\binom{2k+1}{k}   &\cdots  &\binom{2k+1}{2k-3} \\
\hdotsfor4 \\
\hdotsfor4 \\
\binom {2k+1}{2}  &\binom{2k+1}{3}   &\cdots  &\binom{2k+1}{k}
\end{pmatrix}.
$$
By Lemma \ref{det}, the determinant of this matrix is 
$$ 
\frac{\binom{2k+1}{k}\binom{2k+2}{k}\cdots\binom{3k-1}{k}}
{\binom{k}{k} \binom{k+1}{k}\cdots\binom{2k-2}{k}}
$$
which is invertible since the characteristic of the field is $p=3k+1$. Hence
all monomials of degree $3k$ in $A$ and $B$ are elements of the ideal $I$.
\qed \end{proof}

\begin{acknowledgement}
It is a pleasure to thank Mel Hochster for several valuable discussions on
tight closure theory.
\end{acknowledgement}

\end{document}